\documentclass[10pt]{amsart}

\usepackage{color,graphics,graphicx,shortvrb}
\usepackage{epsfig}
\usepackage{newlfont}
\usepackage{caption}

\usepackage{amssymb,amsmath,latexsym}
\usepackage{mathptmx}
\usepackage{subfig}
\usepackage{bibentry}

\newtheorem{thm}{Theorem}[section]

\newtheorem{rem}{Remark}[section]

\newcommand{\R}{\mathbb{R}}
\newcommand{\N}{\mathbb{N}}

\newcommand{\la}{\lambda}
\newcommand{\vf}{\varphi}

\begin{document}

\title[Stabilization of Cascade ODE-Linearized KdV System]{Exponential Stabilization of Cascade ODE-Linearized KdV System by Boundary Dirichlet
 Actuation}

\author{Habib Ayadi}
\address{Universit\'e de Kairouan. Institut Sup\'erieur des Math\'ematiques Appliqu\'ees et de l'Informatique, Avenue Assad Iben Fourat - 3100 Kairouan, Tunisia}
\email{hayadi26@laposte.net}

\begin{abstract}
In this paper, we solve  the problem of exponential stabilization for a class of cascade ODE-PDE system governed by a linear ordinary differential
equation  and a $1-d$ linearized Korteweg-de Vries equation (KdV) posed on a bounded interval.
The control for the entire system acts on the left  boundary with Dirichlet condition of the KdV equation
whereas the KdV acts in the linear ODE by a Dirichlet connection.
 We use the so-called backstepping design in infinite dimension to convert the system under consideration into an exponentially stable
 cascade ODE-PDE system.Then, we use the invertibility of such design  to achieve the exponential stability  for the ODE-PDE cascade system
 under consideration by using Lyapunov analysis.
\end{abstract}
	

\maketitle

\tableofcontents

\thispagestyle{empty}


\section{Introduction }
It is well known  that the Korteweg-de Vries (KdV) equation in bounded domain models
the dynamics of various types of extreme waves in shallow water, more particularly tsunami waves and freak waves.
From theoretical point of view, the KdV controlled equation has some interesting control properties depending on where
the controls are located \cite{Cerp13coro}, \cite{Ros97}. In the past decades,  stabilization of coupled ODE-PDE systems  is widely studied in the literature.
Such systems can be used to model various processes such as road traffic \cite{Goa06}, gas flow pipeline
\cite{Pri15}, power converters connected to transmission lines
\cite{Daaf14}, oil drilling \cite{Sal16} and many other systems.
Many problems of state and output feedback stabilization for  coupled  ODE-Heat has been solved \cite{Tang11}, \cite{Tang11SCL},
\cite{Kris09heat} and ODE- Wave \cite{Kris10JFI}, \cite{Zhoua12},  to cite few.
The problem of controllability of coupled ODE-PDE systems has been discussed in \cite{Weiss09} and \cite{Weiss11}.
Some nonlinear extensions are studied in \cite{Wu13}, \cite{Ahmed15}, \cite{Cai16Kyber} where the non linearity is
assumed to be global Lipshitz, and in \cite{Kris14}, \cite{Hasana16}, \cite{Cai16} for  more general nonlinear ODE.
In this paper, we deal with the state stabilization problem for a cascade ODE-KdV system.
We applied the infinite dimensional backstepping method to build a stabilizing feedback
control for system (\ref{sysdimfini})-(\ref{bound1xx}). The backstepping method was introduced
firstly for finite dimensional control systems governed by ODE \cite{Kris95}. The first extensions to PDE have
appeared in \cite{Coron98} and \cite{Kris00}. Later, in \cite{Liu03} and \cite{Kris04}, the authors have introduced an invertible integral
transformation that transforms the original parabolic PDE into an asymptotically stable one. Recently, the
backstepping method is used to design a feedback control law for coupled PDE-ODE (see the textbook
\cite{Krisbook} and references therein).
As far as we know, problem of stabilization by backstepping design  for such system when the PDE subsystem is governed by the linearized KdV
equation has not yet been tackled in the literature.
This paper is organized as follows. In Section 2, we present the main result of this paper which
is summarized in Theorem \ref{thm}. In section 3, we formulate the backstepping design of the feedback control law. Section 4 is
devoted to prove Theorem \ref{thm} through three steps as follows. Firstly, we prove the invertibility of the transformation given in previous section.
Secondly, we establish the well posedness of system (\ref{sysdimfini})-(\ref{bound1xx}). Finally, we prove
the exponential stability in the sens of the $H$-norm of solutions of (\ref{sysdimfini})-(\ref{bound1xx})  around the origin.

\section{Problem Formulation and Main Result}

Let $l>0$, we consider the following  cascade ODE-PDE system
\begin{align}
\dot{X}(t)&= AX(t) +Bu(l,t),\quad t>0,\label{sysdimfini}\\
u_t(x,t)&= -u_{x}(x,t)-u_{xxx}(x,t),\quad t>0,\quad x\in (0,l) ,\label{kdv}\\
u(0, t) &= U(t),\quad t>0,\label{bound0}\\
u_{x}(l, t) &= 0,\quad t>0,\label{boundlx}\\
u_{xx}(l, t) &= 0,\quad t>0.\label{bound1xxbis}\\
X(0) = X_0, & \, u(x,0) = u_0(x), \quad x \in (0,l).\label{bound1xx}
\end{align}
where $X(t)\in \mathbb R^n$ is the state of the ODE subsystem, $u(x, t)\in \mathbb R$ is the state of the linear KdV subsystem,
$U(t)\in \mathbb R$ is the control input to the entire system acting in the left boundary $x=0$ of the PDE domain  $(0, l)$,
and  $A\in \mathbb{R}^{n\times n}$ , $B \in \mathbb{R}^{n\times 1}$  such that the pair $(A, B)$ is  controllable.
The control objective is to exponentially stabilize the system (\ref{sysdimfini})-(\ref{bound1xx}) around its zero equilibrium.
Along this paper, the Euclidean norm of a vector $X$ in $\mathbb{R}^n$ and the $L^2$-norm of
 a function $u$ in $L^{2}(0,l)$ are denoted by $|X|$ and
  $$\|u\|=\Big(\int_0^lu^2(x)dx\Big)^{\frac{1}{2}},$$
  respectively.
 Let $ H=\mathbb{R}^n\times L^2(0,l)$ the state space of the  system (\ref{sysdimfini})-(\ref{bound1xx}).
 It is obvious that the vector space $H$ equipped  with  its norm
\begin{equation}\label{normH}
\|(X,u)\|_H=\Big(|X|^2+\|u\|^2\Big)^{\frac{1}{2}},
\end{equation}
is a Hilbert space.
The infinite dimensional backstepping design for coupled ODE-PDE system is to seek a continuous and invertible integral transformation
\begin{eqnarray}
  \Omega &:& H\rightarrow H\nonumber\\
  && (X,u)\mapsto (X,w)\label{omega}
\end{eqnarray}
to convert system
(\ref{sysdimfini})-(\ref{bound1xx}) into  the exponentially stable  target system
\begin{align}
\dot{X}(t)&= (A+BK)X(t) +Bw(l,t),\quad t>0,\label{tsysdimfini}\\
w_t(x,t)&= -w_{x}(x,t)-w_{xxx}(x,t)-\lambda w(x,t) ,\quad t>0, \quad x\in (0,l) ,\label{tkdv}\\
w(0, t) &= 0,\quad t>0,\label{tbound0}\\
w_{x}(l, t) &= 0,\quad t>0,\label{tboundlx}\\
w_{xx}(l, t) &= 0,\quad t>0.\label{tbound1xxbis}\\
X(0) = X_0, \, & w(x, 0) = w_0,\quad x \in (0,l),\label{tbound1xx}
\end{align}
where $K\in \mathbb{R}^{1\times n}$ such that $A+BK$ is Hurwitz, and $\lambda$ is an arbitrary positive number.
Define the dense subspace $\Lambda$ of $H$ by   $$\Lambda= \R^n\times \{u\in H^3(0,l)\mid u(0)=u^{'}(l)=u^{''}(l)=0\}.$$
Our main result is stated in the following theorem which asserts that the cascade ODE-PDE system (\ref{sysdimfini})-(\ref{bound1xx})
 is exponentially stable in the sense of the norm $\|\cdot\|_H$.
\begin{thm}\label{thm}
For any initial condition $(X_0, u_0)\in \Lambda$, the closed loop system (\ref{sysdimfini})-(\ref{bound1xx})
with the feedback law (\ref{feedback}) admits a unique classical solution
 $$ (X, u)\in C([0, +\infty) , \Lambda ) \cap C^1([0,+\infty) ; H)$$
and if $(X_0,u_0) \in H$ then the system admits a unique mild solution
$$ (X,u) \in C([0, +\infty) , H).$$
 Moreover, there exists two constants $c_1>0$ and $c_2>0$ such that for all $(X_0,u_0) \in H$ then
 \begin{equation}\label{expst}
  \|(X(t),u(.,t))\|_H \leq c_1 \|(X_0,u_0)\|_He^{-c_2\, t},\,\,\forall t\geq 0.
 \end{equation}
\end{thm}
\section{Control Design}
The transformation $\Omega$  is postulated in the following form :
\begin{align}
X(t)&=X(t),\\
w(x,t)&= u(x,t)-\int_x^l u(y,t)q(x,y)dy-\varphi(x)X(t) ,\label{w}
\end{align}
where the kernel $q(x,y)\in\R$ and  the gain function $\varphi(x)^T\in \R^n$ are to be determined.
From (\ref{w}), the first three derivatives of $w(x,t)$ with respect to $x$  are given by
\begin{equation}
w_x(x,t) = u_x(x,t)+q(x,x)u(x,t)-\int_x^l u(y,t)q_x(x,y)dy-\varphi^{'}(x)X(t), \label{wx}
\end{equation}
\begin{eqnarray}
w_{xx}(x,t)&=&u_{xx}(x,t)+u(x,t)\frac{d}{dx}q(x,x)+q(x,x)u_x(x,t)+q_x(x,x)u(x,t)-\int_x^l u(y,t)q_{xx}(x,y)dy
\nonumber\\
 &&-\varphi^{''}(x)X(t) \label{wxx}\\ and, \;\;\;
w_{xxx}(x,t)&=&u_{xxx}(x,t)+u(x,t)\frac{d^2}{dx^2}q(x,x)+2u_x(x,t)\frac{d}{dx}q(x,x)+q(x,x)u_{xx}(x,t)+u(x,t)\frac{d}{dx}q_x(x,x)
\nonumber\\
            && +q_x(x,x)u_x(x,t)+q_{xx}(x,x)u(x,t)-\int_x^l u(y,t)q_{xxx}(x,y)dy-\varphi^{'''}(x)X(t),
\label{wxxx}
\end{eqnarray}
 respectively.
Furthermore, the derivative of (\ref{w}) with respect to $t$ is given by
\begin{eqnarray}
 w_t(x,t)&=& u_t(x,t)-\int_x^l u_t(y,t)q(x,y)dy-\varphi(x)\dot{X}(t),\nonumber\\
&=& u_t(x,t)+\int_x^l (u_y(y,t)+u_{yyy}(y,t))q(x,y)dy-\varphi(x)(AX(t) +Bu(l,t)).
\end{eqnarray}
Integrating by parts the above identity, we get
\begin{eqnarray}\label{wt}
w_t(x,t)&=&u_t(x,t)-\int_x^l u(y,t)(q_y(x,y)+q_{yyy}(x,y))dy-q(x,x)(u(x,t)+u_{xx}(x,t))\nonumber\\&&
+q_y(x,x)u_x(x,t)-q_{yy}(x,x)u(x,t)+(q(x,l)+q_{yy}(x,l))u(l,t)\nonumber\\&&
-\varphi(x)(AX(t)+Bu(l,t)),
\end{eqnarray}
where the boundary conditions (\ref{boundlx}) and  (\ref{bound1xx}) has been used.
Thus, from (\ref{w}), (\ref{wx}), (\ref{wxxx}) and (\ref{wt}), for all $\la > 0$, the identity
\begin{eqnarray}\label{twila}
 w_t(x,t)+w_x(x,t)+w_{xxx}+\la w(x,t)&=&-\int_x^l u(y,t)(q_{xxx}(x,y)+q_{yyy}(x,y)+q_x(x,y)+q_y(x,y)+\la q(x,y))dy\nonumber\\&&
  +(q(x,l)+q_{yy}(x,l)-\varphi(x)B)u(l,t)\nonumber\\&&
  +u_x(x,t)(q_y(x,x)+q_x(x,x)+2\frac{d}{dx}q(x,x))\nonumber\\&&
  +u(x,t)(\la+q_{xx}(x,y)+\frac{d}{dx}q_x(x,x)+\frac{d^2}{dx^2}q(x,x))\nonumber\\&&
  -(\varphi(x)(A+\la I_n)+\varphi^{''}(x)+\varphi^{'''}(x))X(t),
\end{eqnarray}
 holds. Moreover, setting $x = 0$ in (\ref{w}) and $x=l$ in (\ref{w}), (\ref{wx}) and (\ref{wxx}), it follows
\begin{align}
w(0,t)&=U(t)-\int_0^l q(0,y)u(y,t)dy-\vf(0)X(t),\\
w(l,t)&=u(l,t)-\vf(l)X(t),\\
w_x(l,t)&=q(l,l)u(l,t)-\vf^{'}(l)X(t),\\
w_{xx}(l,t)&=\frac{dq_x}{dx}(l,l)u(l,t)+q_x(l,l)u(l,t)-\vf^{''}(l)X(t).\label{wxxlt}
\end{align}
 Assume that the  gain function $\vf(x)^T$ defined in $[0, l]$ and the kernel $q(x,y)$ defined in the triangle
 \begin{equation}\label{T}
 T=\{(x,y)\mid x\in [0,l],y\in [x,l]\}
 \end{equation}
satisfy
\begin{align}\label{sysfi}
\nonumber \vf^{'''}(x)+\vf^{'}(x)+\vf(x)(A+\la I_n)=&0,\quad x\in [0,l],\\
\nonumber \vf(l)=&K,\\
\nonumber \vf^{'}(l)=&0,\\
\vf^{''}(l)=&0,
\end{align}
and
\begin{align}\label{sysq}
\nonumber q_{xxx}(x,y)+q_{yyy}(x,y)+q_x(x,y)+q_y(x,y)=&-\la q(x,y),\quad (x,y)\in T,\\
\nonumber q(x,l)+q_{yy}(x,l)=&\varphi(x)B,\quad x\in [0,l], \\
\nonumber q(x,x)=&0,\quad x\in [0,l],  \\
q_x(x,x)=&\frac{\la}{3}(l-x),\quad x\in [0,l],
\end{align}
respectively. Then, from  the identities (\ref{twila})-(\ref{wxxlt}), we obtain  the target system (\ref{tsysdimfini})-(\ref{tbound1xx})
 for every solution of
the closed loop (\ref{sysdimfini})-(\ref{bound1xx}) with the feedback law
\begin{equation}\label{feedback}
U(t)= \int_0^l q(0,y)u(y,t)dy+\vf(0)X(t).
\end{equation}
Notice that the existence of the control law (\ref{feedback}) is a direct consequence of the existence of  kernels $q(x,y)$ and $\vf(x)^T$  that satisfy
(\ref{sysq}) and (\ref{sysfi}), respectively. In the following, we compute explicitly the kernel $\vf(x)^T$ and we prove the existence of the kernel
$q(x,y)$.
To begin with, the solution of the ODE ( \ref{sysfi}) is
\begin{equation}\label{fi}
\vf(x)=(K,0,0)e^{(x-l)M}E,
\end{equation}
where $M$ and $E$ are the constant matrices $$M=\begin{pmatrix}
0 & 0 & -(A+\la I_n)\\
I_n & 0 & -I_n\\
0 & I_n & 0
\end{pmatrix}, \;\; E=\begin{pmatrix}
I_n \\
0\\
0
\end{pmatrix}.$$
On the other hand, to prove the existence of the kernel $q(x,y)$, let us make the change of variables
$$s=x+y\quad,\quad t=y-x, $$ and define $$G(s,t):= q(x,y)=q\Big(\frac{s-t}{2},\frac{s+t}{2}\Big).$$
Then, the function $G$ defined in the triangle $T_0=\{(s,t)/t\in[0,l],s\in[t,2l-t]\}$, satisfies
\begin{align}
 (6G_{tts}+2G_{sss}+2G_s)(s,t)=&-\la G(s,t), \quad (s,t)\in T_0,\label{sysG}\\
(G+G_{ss}+2G{st}+G_{tt})(s,2l-s)=&\varphi(s-l)B,\quad s\in [l,2l],\label{s-l}\\
 G(s,0)=&0 ,\quad s\in[0,2l] ,\label{s0}\\
G_t(s,0)=&\frac{\la}{6}(s-2l),\quad s\in[0,2l].\label{Gt0}
\end{align}
To prove the existence of the kernel $q(x, y)$, we adopt the proof done by \cite[Cerpa-Coron]{Cerp13coro} with slight modification.
At this stage, it is conventional to put the previous system in an integral form (see for example \cite{Cerp13coro}).
Using (\ref{s-l}), we rewrite (\ref{sysG}) in variable $(\eta,\xi)$ as follows
\begin{equation}\label{idGss}
(6G_{tts}+6G_{sss}+12G_{sts}+6G_s)(\eta,\xi)=(-\la G+4G_{sss}+12G_{sts}+4G_s)(\eta,\xi).
\end{equation}
Next, we integrate (\ref{idGss}) with respect to the variable $\eta$ over the interval $(s, 2l-\xi)$ and use the boundary condition (\ref{s-l}).
Then, we integrate the obtained identity with respect to $\xi$ over the interval $(0,\tau)$ and use (\ref{s0}) and (\ref{Gt0}). Finally, by integrating the
acquired identity with respect to $\tau$ over the interval $(0,t)$
and use again the boundary condition (\ref{s0}), it follows
\begin{eqnarray}\label{Gintequ}
G(s,t)&=& \int_0^t \int_0^\tau \vf(l-\xi)B d\xi d\tau-\frac{\la t}{6}(2l-s)-2\int_0^t G_s(s,\tau)d\tau-\int_0^t \int_0^\tau (G_{ss}+G)(s,\xi)d\xi d\tau\nonumber\\&&
-\frac{1}{6}\int_0^t \int_0^\tau \int_s^{2l-\xi} (4G_{sss}+4G_s-\la G+12G_{sts})(\eta,\xi) d\eta d\xi d\tau.
\end{eqnarray}
In order to achieve the existence of the function $G$, we use the classical  method of successive approximations \cite{Krisbookcourse}. To this end,  set
\begin{eqnarray}
G^1(s,t)& = & -\frac{\la t}{6}(2l-s),\nonumber\\
G^{n+1}(s,t)& = &\theta(t)-2\int_0^t G^n_s(s,\tau)d\tau-\int_0^t \int_0^\tau (G^n_{ss}+G^n)(s,\xi)d\xi d\tau\nonumber\\
&&-\frac{1}{6}\int_0^t \int_0^\tau \int_s^{2l-\xi} (4G^n_{sss}+4G^n_s-\la G^n+12G^n_{sts})(\eta,\xi) d\eta d\xi d\tau,
\end{eqnarray}
where
\begin{equation}\label{teta}
 \theta(t)= \int_0^t \int_0^\tau \vf(l-\xi)B d\xi d\tau.
\end{equation}
Using an appropriate  calculation, we get
\begin{equation*}
    G^2(s,t)=\theta(t)+\frac{1}{6}\Big(-\la t^2+\Big(\frac{\la l}{9}-\frac{\la^2 l^2}{18}\Big)t^3+\frac{\la }{18}t^4+\frac{\la^2 }{240}t^5
    +\Big(\frac{l\la^2-\la}{18}\Big)st^3-\frac{\la^2 }{72}s^2 t^3\Big).
\end{equation*}
Since $\vf$ is a continuous function on $[0,l]$, from (\ref{teta}), there exists a constant $\rho>0$ such that $|\theta(t)|\leq \rho t^2.$
Keeping this in mind, and using the classical procedure as in \cite{Cerp13coro}, we  get the following inequality
$$|G^k(s,t)|\leq \sum_{i=0}^k\left(\sum_{j=[\frac{3k-1}{2}]-i}^{3k-1-i}a_{ij}t^j\right)s^i,$$
for all  $k\in \N^\ast$, where the positive coefficients $a_{ij}$ have appropriate decay properties so that the series $\sum_{n=1}^{\infty}G^n(s,t)$
is uniformly convergent in $T_0$.  We underline that the uniform convergence  of
the previous series on the domain $T_0$ is due to the fact that the recursive identity is a function of some integral operators.
 Therefore, the series defines a continuous function $$G(s,t)= \sum_{n=1}^{\infty}G^n(s,t), \forall (s,t)\in T_0.$$
 Thus,  the proof of the existence of $G(s,t)$ solution of the integral equation (\ref{Gintequ}) is achieved.
Having proved the existence of $G(s,t)$ in $T_0$, that of $q(x,y)$ in $T$ follows immediately. In consequence, the state feedback controller
(\ref{feedback}) is well defined.

Let's move to the proof of the main Theorem \ref{thm}.
\section{Proof of Theorem {\ref{thm}}}
The proof is divided in three steps. In the first step, we prove the invertibility of the backstepping transformation $\Omega$ defined
by (\ref{omega}). The existence and the global exponential stability of the solution to the closed loop system (\ref{sysdimfini})-(\ref{bound1xx})
with the state feedback (\ref{feedback}) are  established in the second  and the third step, respectively.
\subsection{First step: Invertibility of $\Omega$ }
It is a straightforward that the transformation $\Omega$ is invertible, and
\begin{eqnarray}
\Omega^{-1} &:& H\rightarrow H\\
  \nonumber&& (X,w)\mapsto (X,u)
\end{eqnarray}
has the following form
 \begin{align}
X(t)&=X(t)\\
u(x,t)&= w(x,t)+\int_x^l w(y,t)h(x,y)dy+\psi(x)X(t),\label{u}
\end{align}
where the kernel function $h(x,y)\in\R$  and the gain function $\psi(x)^T\in\R^n$ are to be determined.
As is done in the study of the direct transformation $\Omega$, the derivatives $u_t$, $u_x$, $u_{xxx}$ are computed
and  system (\ref{sysdimfini})-(\ref{bound1xx}) holds if $h(x,y)$  satisfies
\begin{align}\label{sysh}
\nonumber h_{xxx}(x,y)+h_{yyy}(x,y)+h_x(x,y)+h_y(x,y)=&\la h(x, y),\\
\nonumber h(x,l)+h_{yy}(x,l)=-&\psi(x)B, \\
\nonumber h(x,x)=&0,\\
h_x(x,x)=&\frac{\la}{3}(l-x),
\end{align}
in the triangle $T$ defined by (\ref{T}) and  $\psi(x)$ satisfies in $[0,l]$
\begin{align}\label{syspsi}
\nonumber \psi^{'''}(x)+\psi^{''}(x)+\vf(x)(A+BK)=&0,\\
\nonumber \psi(l)=&K,\\
\nonumber \psi^{'}(l)=0,\\
\vf^{''}(l)=&0.
\end{align}
First of all, the solution of the ODE ( \ref{syspsi}) is
\begin{equation}\label{psi}
\psi(x)=(K,0,0)e^{(x-l)N}E,
\end{equation}
where $N$ and $E$ are two constant matrices given by
$$N=\begin{pmatrix}
0 & 0 & -(A+BK)\\
I_n & 0 & -I_n\\
0 & I_n & 0
\end{pmatrix}, \;\;E=\begin{pmatrix}
I_n \\
0\\
0
\end{pmatrix}.$$
Then, in the same way as for the kernel $q(x, y)$, one can easily prove the existence and the continuity of the kernel $h(x, y)$ in $T$. Moreover,
 $\Omega^{-1}$ is continuous operator on the Hilbert space $H$. This achieves
the first step of the proof of Theorem \ref{thm}.
\subsection{Second step: Well Posedness  }
Since  $\Omega : H\rightarrow H $ is a continuous invertible transformation,
 $\Omega^{-1}$  maps a trajectory of (\ref{tsysdimfini})-(\ref{tbound1xx}) into a trajectory of (\ref{sysdimfini})-(\ref{bound1xx})
 with the state feedback controller  (\ref{feedback}). Hence, to prove that (\ref{sysdimfini})-(\ref{bound1xx}) is well posed, it suffices to establish the
 well posedness of the target system  (\ref{tsysdimfini})-(\ref{tbound1xx}). To this reason,
we consider the target KdV subsystem (\ref{tkdv})-(\ref{tbound1xx}) in the Hilbert state space $L^2(0,l)$ and, we define
the unbounded linear operator   $\Gamma: \mathcal{D}(\Gamma)\rightarrow L^2(0,l)$ by
\begin{equation}
\Gamma(w)=-\la w-w^{'}-w^{'''}
\end{equation}
with domain $\mathcal{D}(\Gamma)=\{w\in H^3(0,l)\mid w(0)=w^{'}(l)=w^{''}(l)=0\}$. A simple computation shows that
 the adjoint operator  $\Gamma^{\ast}$ of  $\Gamma$ is
\begin{eqnarray}
\Gamma^{\ast}(w)&=&-\la w+w^{'}+w^{'''},
\end{eqnarray}
with domain $\mathcal D(\Gamma^{\ast})=\{w\in H^3(0,l)\mid w(0)=w^{'}(0)=0,w(l)+w^{''}(l)=0\}$
It is  obvious that the domains $\mathcal D(\Gamma)$ and $ \mathcal D(\Gamma^{\ast})$ are  dense subspaces in $L^2(0,l)$.
Moreover, by integrations by parts, we get
$$ <\Gamma(w), w>=-\frac{1}{2}(w(l)^2+w^{'}(0)^2)-\la \|w\|^2\leq 0,\,\,\forall w\in \mathcal D(\Gamma),$$
$$ <\Gamma^{\ast}(v),v>=-\frac{1}{2}(v(l)^2+v^{'}(l)^2)-\la \|v\|^2\leq 0,\,\,\forall v\in \mathcal D(\Gamma^{\ast}),$$
where $<,>$ stands for the $L^2(0, l)$ standard  inner product. Thus, both operators $\Gamma$ and $\Gamma^{\ast}$ are dissipative.
Therefore, according to the Lumer-Phillips Theorem's, the operator $\Gamma$ generates a $C^0$ semigroup of contraction in $L^2(0, l)$.
Consequently, for all initial condition $w_0\in  \mathcal D(\Gamma)$, the system (\ref{tkdv})-(\ref{tbound1xx}) has a
unique classical solution
$$ w\in C\Big([0,+\infty);  \mathcal D(\Gamma)\Big)\,\cap \, C^1\Big([0,+\infty); L^2(0,l)\Big).$$
Furthermore, by Duhamel formula's, for all initial condition $X_0\in \R^n$, the ODE (\ref{tsysdimfini}) has a unique global solution
$$X(t)=  e^{t(A+BK)}X_0+\int_0^t e^{(t-\tau)(A+BK)}Bw(l,\tau)d\tau.$$
We conclude that the target system  (\ref{tsysdimfini})-(\ref{tbound1xx}) is well posed in $H$, and
for all initial condition $(X_0, w_0)\in \Lambda$, system (\ref{tsysdimfini})-(\ref{tbound1xx}) has a unique classical solution
$$ (X, w)\in C\big([0,+\infty); \Lambda\big)\,\cap \, C^1\big([0,+\infty); H \big),$$
and for $(X_0, w_0)\in H$, the system (\ref{tsysdimfini})-(\ref{tbound1xx}) has a unique mild solution
$$ (X, w)\in C\big([0,+\infty); H\big).$$
\subsection{Third step: Exponential Stability }
Consider the Lyapunov function candidate
\begin{equation}\label{Vly}
V(t)= X(t)^TPX(t)+\frac{\mu}{2}\|w(.,t)\|_{L^2(0,l)}^2,
\end{equation}
where $\mu>0$ is a  constant to be designed later and the positive definite matrix $P=P^T>0$ is the solution of the Lyapunov equation
$$ P(A+BK)+(A+BK)^T P = -Q,$$
for some positive definite matrix $Q=Q^T>0$.  From (\ref{Vly}), it can be obtained that for all $t\geq 0$,
\begin{equation}\label{Venc}
 \alpha_1\|(X(t), w(.,t))\|_H^2\leq V(t)\leq \alpha_2\|(X(t), w(.,t))\|_H^2,
\end{equation}
where $$\alpha_1=\min\big(\la_{min}(P),\frac{\mu}{2})\,\,\,\,\, and \,\,\,\,\,\alpha_2=\max\big(\la_{max}(P),\frac{\mu}{2}\big).$$
The derivative of $V$ along the solutions of (\ref{tsysdimfini})-(\ref{tbound1xx}) is given by
$$\dot{V}(t)=-X^T(t)QX(t)+2X(t)^T PBw(l,t)-\mu w(l,t)^2-\mu\la \|w(.,t)\|^2.$$
By Young's inequality, we get
$$2X(t)^T P B w(l,t)\leq \frac{\la_{min}(Q)}{2}|X(t)|^2+\frac{2}{\la_{min}(Q)}|PB|^2 w(l,t)^2.$$
Thus, $$ \dot{V}(t)\leq -\frac{\la_{min}(Q)}{2}|X(t)|^2-\Big(\mu-\frac{2}{\la_{min}(Q)}|PB|^2\Big)w(l,t)^2-\mu\la \|w(.,t)\|^2.$$
Now, using (\ref{Venc}) and  choose
\begin{equation}\label{mu}
\mu > Max\Big(\frac{2}{\la_{min}(Q)}|PB|^2, 2\la_{max}(P)\Big).
\end{equation}
Then, the inequality
$$\dot{V}(t)\leq -\delta V(t),$$
holds for all $ t\geq 0$, where $\delta= \dfrac{min(\la_{min}(Q),4
\la)}{\mu}.$
Therefore,
\begin{equation}\label{Vexp}
V(t)\leq V(0)e^{-\delta t}, \forall t\geq 0.
\end{equation}
Using (\ref{Venc}) and (\ref{Vexp}), by tacking $\alpha=\frac{\alpha_2} {\alpha_1}$, it follows that
\begin{equation}\label{wexp}
\|(X(t), w(.,t))\|_H^2\leq \alpha \|(X_0, w_0)\|_H^2e^{-\delta t},\,\,\, \forall t\geq 0.
\end{equation}
Recall that the transformations $\Omega$ and $\Omega^{-1}$ are linear and continuous, then there exist two positive constants
$d_1$ and  $d_2$ such that
\begin{align}
\|(X,w)\|_H=\|\Omega(X,u)\|_H\leq &\,\, d_1 \|(X,u)\|_H,\label{cu}\\
\|(X,u)\|_H=\|\Omega^{-1}(X,w)\|_H\leq & \,\,d_2 \|(X,w)\|_H.\label{cw}
\end{align}
Hence, for all initial condition $\big(X_0,u_0\big)\in H$, for all $t\geq 0$, we obtain from (\ref{wexp}), (\ref{cu}) and (\ref{cw})
\begin{equation}
\|(X(t),u(.,t))\|_H\leq c_1 \|(X_0, u_0)\|_He^{- c_2 t},
\end{equation}
where $c_1=d_1d_2\sqrt{\alpha}$ and  $c_2=\frac{\delta}{2}.$
Thus, the proof of Theorem \ref{thm} is complete.
\begin{rem}
Since the positive definite matrix $Q$ and the positive parameter $\lambda$ are arbitrary, the decay rate $c_2$ in (\ref{expst})
is arbitrary. Thus, the system (\ref{sysdimfini})-(\ref{bound1xx}) is rapidly exponentially stable.
\end{rem}

\end{document}